\documentclass[preprint,12pt,authoryear]{elsarticle}

\usepackage{amssymb}
\usepackage{amsmath}

\journal{arXiv}

\begin{document}

\begin{frontmatter}



\title{From Reflecting Brownian Motion to Reflected Stochastic Differential Equations: A Systematic Survey and Complementary Study}


\author[1]{Yunwen Wang}
\author[2]{Jinfeng Li}

\address[1]{Department of Statistics, University of Warwick, Coventry, CV4 7AL, UK}
\address[2]{Department of Electrical and Electronic Engineering, Imperial College London, London SW7 2AZ, UK}

\begin{abstract}
\noindent This work contributes a systematic survey and complementary insights of reflecting Brownian motion and its properties. Extension of the Skorohod problem’s solution to more general cases is investigated, based on which a discussion is further conducted on the existence of solutions for a few particular kinds of stochastic differential equations with a reflected boundary. It is proved that the multidimensional version of the Skorohod equation can be solved under the assumption of a convex domain (D).
\end{abstract}




\end{frontmatter}


\section{Introduction}
\noindent Brownian motion and reflecting Brownian motion have a set of properties that hold great potential for game-changing applications, ranging from the mathematical study of queuing models with heavy traffic (\cite{reiman1984open}), to statistical physics (\cite{lang1995effective}), and the recent advances in statistical mechanics (\cite{grebenkov2019probability}). Underpinned by the It$\hat{o}$'s calculus (\cite{malliaris1983ito}), the Brownian motion inspired models arguably lay a foundation for quantitative finance, as evidenced notably in the stock market forecasting (\cite{guo2019novel}), the valuations of options based on the Black–Scholes–

\noindent Merton model (\cite{calin2012introduction}) in lattice-based computation or Monte Carlo based derivative investment instruments (\cite{zhang2020value}) with multiple sources of uncertainty.

\noindent Fundamentally, the reflecting Brownian motion can be characterised by the Skorohod problem into solving a stochastic differential equation (SDE) with a reflecting boundary condition. Over decades, there has been a continuous research campaign to attempt solutions of the Skorohod stochastic differential equation with reflecting (\cite{saisho1987stochastic}), semi-reflecting (\cite{kanagawa2009numerical}), or absorption (\cite{berestycki2014critical}) boundary conditions. By way of illustration, (\cite{d2012markov}) reports a reflected Markov-modulated Brownian motion with a two-sided reflection, which generalises the reflected Brownian motion to the Markov modulated case. More recently, Monte Carlo simulation was employed by (\cite{malsagov2019approximations}) for closed-form approximations of the mean and variance of fractional Brownian motion reflected at level 0, the problem of which explicit expressions and numerical methods struggle to address. Apart from the discretization error reported by (\cite{asmussen1995discretization}), there remain massive technological gaps that challenge the conventional statistical thinking in tailoring the reflecting Brownian motion and its properties for real-life emerging applications, such as the contact tracing (\cite{li2020global}) for the coronavirus disease (COVID-19), the statistical clutter modelling and phased array signal processing for 5G communications (\cite{li2020low}) and beyond.

\noindent This paper is organised into the following sections. First, the elements of Brownian Motion and the stochastic integral are elaborated in sections 2 and 3, respectively. Section 4 presents the survey and insights into the reflecting Brownian motion, followed by the in-depth discussion of reflected stochastic differential equations in section 5.

\section{Brownian Motion}

\noindent In this section, we will introduce the definition of Brownian motion and reflecting brownian motion, which are the most important stochastic process that is widely used in applications.

\subsection{Brownian motion}

\noindent\textbf{Definition 2.1.1(Brownian Motion)}: A stochastic process X=($X_t$,t$\geq$0) is called a d-dimensional \textbf{Brownian motion}(or \textbf{Wiener process}) with the initial probability law $\mu$, if

(i) $X_t$ is continuous in t almost surely and $X_0$ has the distribution law $\mu$ on $\mathbb{R}^d$;

(ii) $X_t$ has independent increments, that is for any 0$\leq t_0< t_1<...<t_n$, the random variables,

\begin{center}
    $X_{t_n}-X_{t_{n-1}},X_{t_{n-1}}-X_{t_{n-2}},...,X_{t_2}-X_{t_1}$
\end{center}

\noindent are independent with each other;

(iii) for all 0$\leq s<t$,

\begin{center}
    $X_t-X_s\sim$N(0,t-s)
\end{center}

\noindent Then, for every 0$<t_1<...<t_m$ and $A_i\in \mathcal{B}(\mathbb{R}^d)$, i=1,2,...,m, we have the probability as

\begin{center}
    P($X_{t_1}\in A_1,X_{t_2}\in A_2,...,A_{t_m}\in A_m$)
        
    =$\int_{\mathbb{R}^d}\mu(dx)\int_{A_1}p(t_1,x_1-x)dx_1\int_{A_2}p(t_2-t_1,x_2-x_1)dx_2$
        
    ...$\int_{A_m}p(t_m-t_{m-1},x_m-x_{m-1})dx_m$
\end{center}

\noindent where, p(t,x),$t>0$, x$\in \mathbb{R}^d$, is defined by

\begin{center}
    p(t,x)=(2$\pi t)^{-d/2}exp[-|x|^2/2t]$,
\end{center}

\noindent being the probability distribution function of a d-dimensional Gaussian distribution.

\noindent Then we are going to introduce a new concept called Brownian local time which will be used in the upcoming proof.

\noindent Let X=($X_t$) be a one-dimensional Brownian motion defined on the probability space ($\Omega$,F,P).\\
\\

\noindent\textbf{*Definition 2.1.2(Brownian local time)}\cite{ikeda2014stochastic}:By the \textbf{local time} or the \textbf{sojourn time density} of X we mean family of non-negative random variables {$\phi(t,x,\omega),t\in[0,\infty),x\in R^1$} such that, with probability one, the following holds:

(i) (t,x) $\mapsto \phi$(t,s) is continuous,

(ii) for every Borel subset A of $R^1$ and t$\geq$0

$\int^t_0I_A(X_s)ds$=2$\int_A\phi(t,s)dx$.

\noindent It is clear that if such a family \{$\phi$(t,x)\} exists, then it is unique and is given by

$\phi(t,x)=lim_{\epsilon\downarrow 0}\frac{1}{4\epsilon}\int^t_0 I_{(x-\epsilon,x+\epsilon)}(X_s)ds$.

\subsection{Skorohod problem and Skorohod equation}

\noindent Before moving on to reflecting Brownian motion, we will introduce a concept called Skorohod problem. It is useful when characterising the reflecting Brownian motion.

\noindent In probability theory, the Skorokhod problem is the problem of solving a stochastic differential equation with a reflecting boundary condition. Let us firstly see the one dimensional version.

\noindent We set \textbf{$W_0^1$} = \{f$\in$$\mathbb{C}$([0,$\infty)$$\longrightarrow$$\mathbb{R}$); f(0)=0\} and \textbf{$C^+$} = \{f$\in\mathbb{C}([0,\infty$)$\longrightarrow$$\mathbb{R}$); f(t)$\geq$0 for all t$\geq$0\}.

\noindent \textbf{Lemma 2.2.1} \cite{ikeda2014stochastic}: Given f$\in$\textbf{$W_0^1$} and x$\in$\textbf{$\mathbb{R}^+$}, there exist unique g$\in$\textbf{$C^+$} and h$\in$\textbf{$C^+$} such that

(i) g(t)=x+f(t)+h(t),

(ii)h(0)=0 and t$\mapsto$h(t) is increasing,

(iii) $\int^t_0$\textbf{$I_{0}$}(g(s))dh(s)=h(t), i.e., h(t) increases only on the set of t when g(t)=0.

\noindent Proof:

Set

\begin{center}
    g(t)=x+f(t)-$min_{0\leq s\leq t}$\{(x+f(s))$\land$0\},
    
    h(t)=-$min_{0\leq s\leq t}$\{(x+f(s))$\land$0\}.
\end{center}

\noindent Then we will prove that g(t) and h(t) satisfy the above conditions (i), (ii) and (iii). According to the definition of g(t) and h(t), (i) apparently holds. So we only need to check (ii) and (iii).

For (ii), h(0)=-\{(x+f(0)$\land$0\}.

For (ii), h(0)=-\{(x+f(0)$\land$0\}.

x\textgreater0 and f(0)=0 since x $\in R^+$, f $\in W_0'$ $\implies$ x+f(0)\textgreater0

So h(0) = -0 = 0.

\noindent Assume $t_1 < t_2$, $t_1, t_2 \in [0,\infty)$. Then $min_{0\leq s\leq t_1}$\{(x+f(s)) $\land$0\} $\leq$ $min_{0\leq s\leq t_2}$\{(x+f(s)) $\land$0\} because [0, $t_1$] $\subseteq$ [0, $t_2$]. So h($t_1$)=$min_{0\leq s\leq t_1}$\{(x+f(s)) $\land$0\} $\leq$ $min_{0\leq s\leq t_2}$\{(x+f(s)) $\land$0\}=h($t_2$), i.e. t $\mapsto$h(t) is increasing.

\noindent Hence (ii) holds.

\noindent To prove (iii), we need to prove that the set I, which is the union of the intervals on which h(t) increases, is included by $I_0$(g(t)).

\noindent Assume that, $min_{0\leq s\leq t}$\{x+f(s)\} decreases on interval I', i.e. $min_{0\leq s\leq b}$\{x+f(s)\} $\leq$$min_{0\leq s\leq a}$\{x+f(s)\} for all b\textgreater a, where a,b $\in$I'. When b$\longrightarrow$a,

\begin{center}
    -(x+f(b))=-$min_{0\leq s\leq b}$\{x+f(s)\}=h(b),
\end{center}

and 

\begin{center}
    -$min_{0\leq s\leq a}$\{x+f(s)\}=h(a)= $lim_{b\longrightarrow a}$h(b)= $lim_{b\longrightarrow a}$-(x+f(b))=-(x+f(a)),
\end{center}

\noindent since h(t) and f(t) are both continuous.

\noindent Therefore, g(a)=0, g(b)=0. Then because a,b $\in$I' are arbitrary, g(t)=0 on I', i.e. I $\subseteq I_0(g(t))$. Hence (iii) holds.

\noindent We shall prove the uniqueness. Suppose $g_1$(t) and $h_1$(t) $\in C_+$ also satisfy the condition (i), (ii) and (iii). Then

\begin{center}
    g(t)- $g_1$(t)=h(t)- $h_1$(t)
    for all t $\geq$0.
\end{center}

\noindent If there exists $t_1>$0 such that g($t_1$)-$g_1(t_1)>$0, we set $t_2$= max\{ t$<t_1$; g(t)-$g_1$(t)=0\}. Then g(t)$>g_1$(t)$\geq$0 for all t $\in(t_2,t_1]$ and hence, by (iii), h($t_1$)-h($t_2$)=0. Since $h_1$(t) is increasing, we have

\begin{center}
    0$<g(t_1$)-$g_1(t_1)$=h($t_1$)-$h_1(t_1)\leq$h($t_2$)-$h_1(t_2)$=g($t_2$)-$g_1(t_2)$=0.
\end{center}

\noindent This is contradiction. Therefore g(t)$\leq g_1$(t) for all t $\geq$0. By symmetry, g(t)$\geq g_1$(t) for all t $\geq$0. Hence g(t)$\equiv g_1$(t) and so h(t)$\equiv h_1$(t).

\subsection{Reflecting Brownian motion}

\noindent From the name, we can easily see that this Brownian motion has bounds and gets reflected when it goes beyond the bounds.

\noindent\textbf{Definition 2.3.1(Reflecting Brownian motion)} \cite{ikeda2014stochastic}: Let X=($X_t$) be a one-dimensional Brownian motion and let $X^+$=($X^+_t$) be a continuous stochastic process on [0,$\infty$) defined by

\begin{center}
    $X^+_t=|X_t|$.
\end{center}

\noindent Then we can easily infer that the event \{$X^+_{t_2}=y|X_{t_1}=x$\} is the union of \{$X_{t_2}=y|X_{t_1}=x$\} and \{$X_{t_2}=-y|X_{t_1}=x$\}, i.e.

\begin{center}
    \{$X^+_{t_2}=y|X_{t_1}=x$\}=\{$X_{t_2}=y|X_{t_1}=x$\}$\cup$\{$X_{t_2}=-y|X_{t_1}=x$\}
    
    =\{$X_{t_2-t_1}=y-x$\}$\cup$\{$X_{t_2-t_1}=-y-x$\}.
\end{center}

\noindent What's more, since $X_{t-s}\sim$N(0,t-s), we have that for 0$<t_1<t_2<...<t_n$, $A_t\in\mathcal{B}([0,\infty$)),

\begin{center}
    P[$X^+_{t_1}\in A_1$,$X^+_{t_2}\in A_2$,...,$X^+_{t_n}\in A_n$]
    
    =$\int_{[0,\infty)}\mu(dx)\int_{A_1}p^+(t_1,x,x_1)dx_1\int_{A_2}p^+(t_2-t_1,x_1,x_2)dx_2$
    
    ...$\int_{A_1}p^+(t_n-t_{n-1},x_{n-1},x_n)dx_n$
\end{center}

\noindent where

\begin{center}
    $p^+(t,x,y)=\frac{1}{\sqrt{2\pi t}} (exp[-\frac{(x-y)^2}{2t}]+exp[-\frac{(x+y)^2}{2t}])$
\end{center}

\noindent and $\mu^+$ is the probability law of $X_0^+=|X^+_0|$. The process $X^+$ is called the one-dimensional reflecting Brownian motion. An illustrated one-dimensional Brownian motion and the associated reflected path are plotted in Figure 1.\\
\\

\begin{figure}[h]
    \centering
    \includegraphics[scale=0.5]{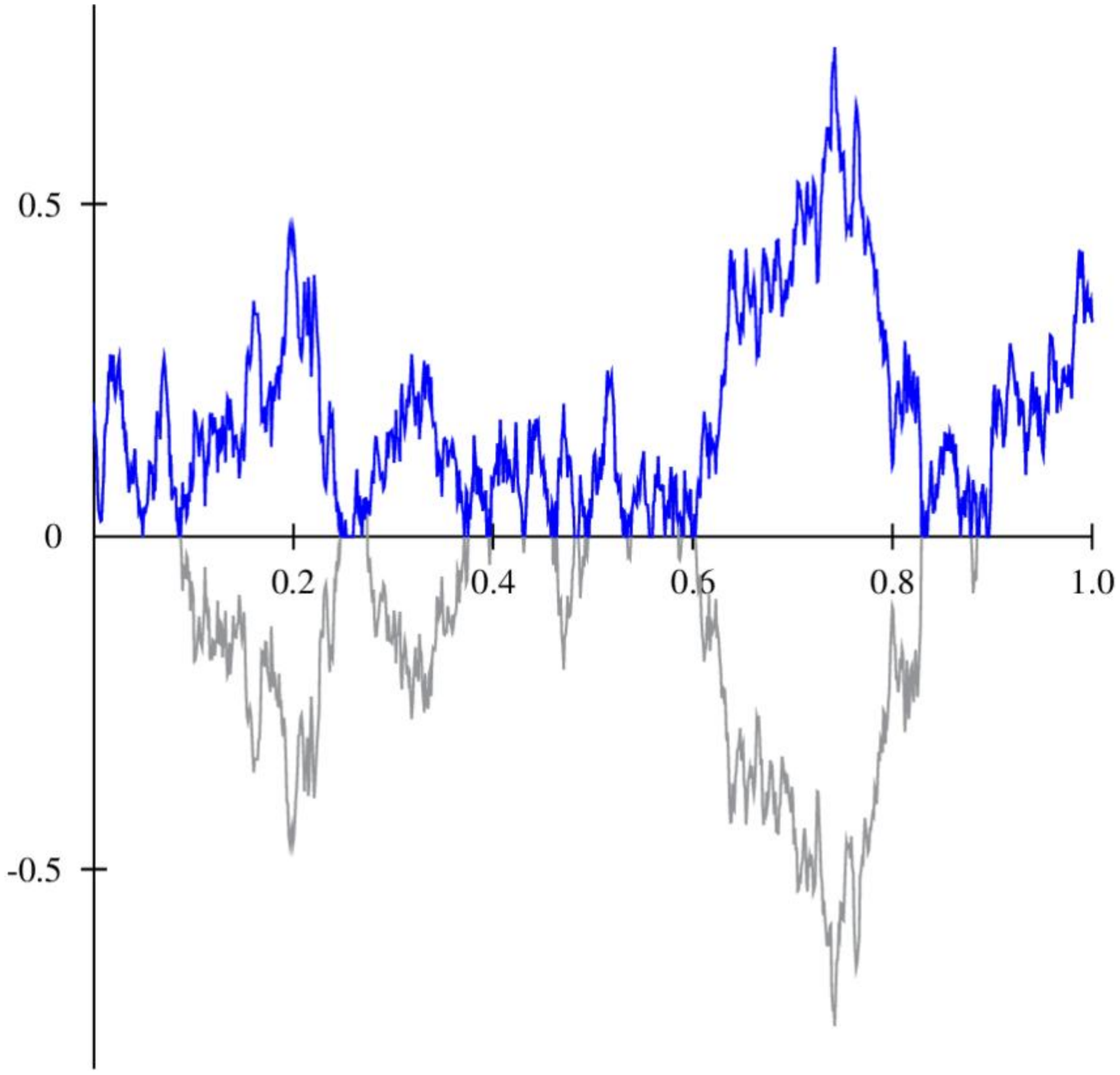}
    \caption{A path of a Brownian motion (grey) and the associated reflected path (blue)}
    \label{fig:my_label}
\end{figure}

\noindent There are different ways to characterise reflecting Brownian motion. We will present a characterisation due to Skorohod problem.

\noindent \textbf{Theorem 2.3.2} \cite{ikeda2014stochastic}: Let \{X(t), B(t),$\phi$(t)\} be a system of real continuous stochastic processes defined on a probability space such that B(t) is a one-dimensional Brownian motion with B(0)=0, X(0) and process {B(t)} are independent and with probability one the following holds:

(i) X(t)$\geq$0 for all t $\geq$0 and $\phi$(t) is increasing with $\phi$(0)=0 such that

\begin{center}
    $\int_0^t I_{0}$(X(s))d $\phi$(s) = $\phi$(t);
\end{center}

(ii)
\begin{equation}
    X(t) = X(0)+B(t)+\phi(t).
\end{equation}

\noindent Then X=X(t) is a reflecting Brownian motion on [0,$\infty$).

\noindent Equation (2.1) is called the \textbf{Skorohod equation} which will be introduced in the next section.

\noindent Proof: By Lemma 2.2.1, X = X(t) and $\phi =\phi$(t) are uniquely determined by X(0) and B = B(t): X = X(0)+B(t)-$min_{0\leq s\leq t}$\{(X(0)+B(s))$\land$0\} and $\phi$= -$min_{0\leq s\leq t}$\{(X(0)+B(s))$\land$0\}. To prove this theorem we only need to show that if $x_t$ is a one-dimensional Brownian motion, then X(t) = $|x_t|$ satisfies, with some processes B(t) and $\phi$(t), the above properties. Let $g_n$(x) be a non-negative continuous function on $R^1$ with support in (0,1/n), i.e. supp($g_n$)=\{x $\in [0,\infty) |g_n$(t)$\neq$0\}=(0,1/n), such that $\int_0^{\infty}g_n$(x)dx=1. Set

\begin{center}
    $u_n$(x) = $\int_0^{|x|}$dy $\int_0^y g_n$(z)dz.
\end{center}

\noindent Then we can see that $u_n \in C^2$(R), $u_n$'(x) =$\frac{d|x|}{dx}\int_0^{|x|}g_n(z)dz$. And then $|u_n'|$ $\leq |\frac{d|x|}{dx}||\int_0^{\infty}g_n(z)dz| \leq$1, $u_n$(0)=0, $u_n$(x) $\uparrow |x|$. Then because $lim_{n\longrightarrow \infty} \frac{1}{n}$=0 and non-negative function $g_n$ has support in (0,1/n),

\begin{center}
    $lim_{n\longrightarrow\infty}\int_0^{|x|}g_n$(z)dz = $lim_{n\longrightarrow\infty}(\int_0^{1/n}g_n$(z)dz+$\int_{1/n}^{|x|}g_n$(z)dz) = 1+0=1 , x$\neq$0
\end{center}

\noindent so $u_n$'(x) = $\frac{d|x|}{dx}$ , for x $\neq$0. Hence,

\begin{equation}
u_n'(x)= sgn(x)= \left\{
        \begin{array}{cc}
             1, & x>0 \\
             0, & x=0 \\
             -1, & x<0
        \end{array}
\right.
\end{equation}

\noindent By It$\hat{o}$'s formula (will be introduced in the following section),

\begin{center}
    $u_n(x_t)$-$u_n(x_0)$=$\int_0^t u_n'(x_s)dx_s$+$\frac{1}{2}\int_0^t u_n''(x_s)$ds
    
    =$\int_0^t u_n'(x_s)dx_s$+$\int_{-\infty}^0 g_n(-y)\phi(t,y)dy$+$\int_0^{\infty}g_n(y)\phi(t,y)dy$,
\end{center}

\noindent where $\phi$(t,y) is the local time of $x_t$. Letting n $\longrightarrow \infty$, we have

\begin{center}
    X(t)-X(0) = $\int_0^t$sgn($x_s)dx_s$+2$\phi$(t,0).
\end{center}

\noindent Set

\begin{center}
    B(t) = $\int_0^t$sgn($x_s$)$dx_s$ and $\phi$(t) = 2 $\phi$(t,0).
\end{center}

\noindent Then, since $<B>_t$=t, B(t) is an (\textbf{$F_t$})-Brownian motion, where (\textbf{$F_t$})=(\textbf{$F_t^x$}) is the proper reference family for $x_t$. Since X(0) is \textbf{$F_0$}-measurable, X(0) and \{B(t)\} are independent. Since

\begin{center}
    $\phi$(t) = $lim_{\epsilon\downarrow 0}\frac{1}{2\epsilon}\int_0^t I_{[0,\epsilon)}$(X(s))ds,
\end{center}

\noindent it is clear that

\begin{center}
    $\int_0^t I_0$(X(s))d$\phi$(s) = $\phi$(t).
\end{center}

\noindent Therefore \{X(t), B(t), $\phi$(t)\} satisfies all conditions in Theorem 2.1. Thus X-(X(t)) and $\phi$=($\phi$(t)) are characterized as X = X(0)+B(t)-$min_{0\leq s\leq t}$\{(X(0)+B(s))$\land$0\} and $\phi$ = -$min_{0\leq s\leq t}$\{(X(0)+B(s))$\land$0\}.\\
\\
\noindent\textbf{Theorem 2.3.3} \cite{ikeda2014stochastic}: The local time \{$\phi$(t,x)\} of X exists.

\noindent\textbf{Proof}: We will prove this theorem by using stochastic calculus. Let ($F_t)=(F_t^X)$ be the proper reference family of X. Then X is an ($F_t$)-Brownian motion and $X_t-X_0$ belong s to space M. Let $g_n$(x) be a continuous function on $R^1$ such that its support is contained in (-1/n+a,1/n+a), $g_n(x)\geq0$, $g_n$(a+x)=$g_n$(a-x) and

$\int_{-\infty}^\infty g_n(x)dx$=1.

\noindent Set

$u_n(x)=\int_{-\infty}^x dy\int_{-\infty}^y g_n(z)dz.$

\noindent By It$\hat{o}$'s formula,

$u_n(X_t)-u_n(X_0)=\int_0^tu_n'(X_s)dX_s+\frac{1}{2}\int_0^tu_n''(X_s)ds$

\noindent and if the local time \{$\phi(t,x)$\} does exist, then

$\frac{1}{2}\int_0^tu_n''(X_s)ds=\frac{1}{2}\int_0^tg_n(X_s)ds=\int_{-\infty}^\infty g_n(y)\phi(t,y)dy\longrightarrow\phi(t,a)$ as n$\longrightarrow\infty$.

\noindent Also, it is clear that

$u_n(x)\longrightarrow(x-a)^+$,
$u_n'(x)\longrightarrow
\left\{
        \begin{array}{cc}
             1, & x>a \\
             \frac{1}{2}, & x=a \\
             0, & x<a
        \end{array}
\right.$

\noindent Hence, $\phi(t,a)$ should be given as

\begin{equation}
    \phi(t,a)=(X_t-a)+-(X_0-a)^+-\int_0^tI_{(a,\infty)}(X_s)dX_s.
\end{equation}
\\
\\
\noindent However, the theorems stated above raise a new problem, i.e. the integration of stochastic process. Note that a continuous stochastic process can be nowhere differentiable, and the necessary condition of ordinary integration is not satisfied, hence the need to find a new integral. A way around the obstacle was found by It$\hat{o}$ in the 1940s which will be introduced in the next section.

\section{Stochastic Integral}

\subsection{It$\hat{o}$ Integral}

\noindent We only give the relative definition and theorem of stochastic integral with respect to Brownian motion
$\int_0^T$f(t)d$B_t$
which is also called It$\hat{o}$ stochastic integral.

\noindent The It$\hat{o}$ integral is a random variable since $B_t$ and the integrand f(t)(to be precise, f(t,$\omega$)) are random. In order to guarantee the regularity of the integral, we will give some restrictions on f(t).

\noindent\textbf{Definition 3.1.1}: Denote $\mathcal{L}^2$ to be the set of random variables X, in which $E|X|^2<\infty$.

\noindent\textbf{Definition 3.1.2($\mathcal{M}^2$ Stochastic Process)}: Denote $\mathcal{M}^2$ to be the class of stochastic processes f(t), t$\geq$0, such that

$E(\int_0^T|f(t)|^2dt)<\infty$.

\noindent Let $\mathcal{M}^2_T$ be the class of stochastic processes f(t) such that f(t)$\in\mathcal{M}^2_T$ for any T>0.

\noindent\textbf{Definition 3.1.3($\mathcal{L}^2$ and $\mathcal{M}^2_T$ Norm)}: For a random variable X the $\mathcal{L}^2$ norm is $||X||_{\mathcal{L}^2}=\sqrt{E(X^2)}$.

\noindent For a stochastic process f$\equiv$f(t) the $\mathcal{M}^2_T$ norm is $||f||_{\mathcal{M}^2_T}=\sqrt{E(\int^T_0|f(t)|^2dt)}$.

\noindent\textbf{Definition 3.1.4($\mathcal{L}^2$ and $\mathcal{M}^2_T$ Convergence)}: A sequence of random variables \{$X_n$\} converges in $\mathcal{L}^2$ to X if

$||X_n-X||_{\mathcal{L}^2}=E(|X_n-X|^2)\longrightarrow0$.

\noindent A sequence of random functions/stochastic processes \{$f_n(t)$\} converges in $\mathcal{M}^2_T$ to f if 

$||f_n-f||_{\mathcal{M}^2_T}=E(\int_0^T|f_n(t)-f(t)|^2 )\longrightarrow$0.

\noindent Similarly, \{$f_n(t)$\} converges in $\mathcal{M}^2$ to f if \{$f_n(t)$\} converged to f in $\mathcal{M}_T^2$ for all T.

\noindent Based on definition 3.1.1,3.1.2,3.1.3 and 3.1.4, It$\hat{o}$ integral can be defined on $\mathcal{M}^2$.

\noindent\textbf{Definition 3.1.5(It$\hat{o} Integral$} For any T\textgreater0 and any stochastic process f$\in\mathcal{M}^2$, the stochastic integral of f on [0,T] is defined by

$I_T(f)=\int_0^Tf(t)dB_t$.

\noindent\textbf{Theorem 3.1.6(Existence and Uniqueness of It$\hat{o}$ Integral)}: Suppose that a function f$\in\mathcal{M}^2$ satisfies the following assumptions

(i) f(t) is almost surely continuous, i.e. P($lim_{\epsilon\longrightarrow0}|f(t+\epsilon)-f(t)|=0$)=1;

(ii) f(t) is adapted to the filtration \{$\mathcal{F}_t$\}, where $\mathcal{F}_t=\sigma(\{B_s,s<t\}$).

\noindent Then, for any T\textgreater0, the It$\hat{o}$ integral

$I_T(f)=\int_0^Tf(t)dB_t$

\noindent exists and is unique almost everywhere.\\
\\
\noindent\textbf{Example 3.1.1}: To show the existence of $\int^T_0B_t^2dB_t$, we need to show that the $B_t^2$ belongs to $\mathcal{M}^2$. Since for all T

$E(\int^T_0|B_t|^4dB_t)=\int^T_0E(|B_t|^4)dB_t=\int^T_03t^2dt<\infty$.

\noindent So $B_t^2\in\mathcal{M}^2$. Also it is easy to verify that it satisfies (i) and (ii) of Theorem 3.1.6. Hence the It$\hat{o}$'s Integral exists.\\
\\
\noindent\textbf{Theorem 3.1.7}: The following properties holds for any f, g$\in\mathcal{M}^2$, any $\alpha,\beta\in\mathbf{R}$ and any 0$\leq s<$t:

(i)Linearity:

\begin{equation}
    \int_0^t(\alpha f(r)+\beta g(r))dB_r)=\alpha\int_0^tf(r)dB_r+\beta\int_0^tg(r)dB_r;
\end{equation}

(ii)Isometry:

\begin{equation}
    E(|\int_0^tf(r)dB_r|^2)=E(\int_0^t|f(r)|dr);
\end{equation}

(iii)Martingale Property:

$E(\int_0^tf(r)dB_r|\mathcal{F}_s)=\int_0^sf(r)dB_r$.
In particular, $E(\int_0^tf(r)dB_r)$=0.

\subsection{It$\hat{o}$'s Lemma}

\noindent\textbf{Definition 3.2.1(It$\hat{o}'s$ Lemma)}: Suppose that F(t,x) is a real valued function with continuous partial derivatives $F_t$(t,x), $F_x$(t,x) and $F_{xx}$(t,x) for all t$\geq$0 and x$\in\mathbf{R}$. Assume also that the process $F_x(t,B_t)\in\mathcal{M}^2$. Then F(t,$B_t$) satisfies

\begin{equation}
    F(T,B_t)-F(0,B_0)
    =\int_0^T[F_t(t,B_t)+\frac{1}{2}F_{xx}(t,B_t)]dt+\int_0^TF_x(t,B_t)dB_t,
\end{equation}
a.s.

\noindent In differential notation, (3.3) can be written as

\begin{equation}
    dF(t,B_t)=[F_t(t,B_t)+\frac{1}{2}F_{xx}(t,B_t)]dt+F_x(t,B_t)dB_t.
\end{equation}

\noindent\textbf{Proof}: We first prove only the case where F, $F_x, F_{xx}$ are all bounded by some C\textgreater0. Consider a partition of [0,T], 0=$t_0^n<t_1^n<...<t_n^n$=T, where $t_i^n=\frac{iT}{n}$. Denote $B_{t_i^n}$ by $B_i^n$; the increments $B_{i+1}^n-B_i^n$ by $\Delta_i^nB$; and $t_{i+1}^n-t_i^n$ by $\Delta_i^nt$. Using Taylor's expansion, there is a point $\hat{B}_i^n$ in each interval [$B_i^n,B_{i+1}^n$] and a point $\hat{t}_i^n$ in each interval [$t_i^n,t_{i+1}^n$] such that

$F(T,B_T)-F(0,B_0)$

$=\Sigma_{i=0}^{n-1}[F(t_{i+1}^n,B_{i+1}^n)-F(t_i^n,B_i^n)]$

$=\Sigma_{i=0}^{n-1}[F(t_{i+1}^n,B_{i+1}^n)-F(t_i^n,B_{i+1}^n)]+\Sigma_{i=0}^{n-1}[F(t_i^n,B_{i+1}^n)-F(t_i^n,B_i^n)]$,

and then by Taylor's expansion

$=\Sigma_{i=0}^{n-1}F_t(\hat{t}_i^n,B_{i+1}^n)\Delta_i^nt+\Sigma_{i=0}^{n-1}F_x(t_i^n,B_i^n)\Delta_i^nB+\frac{1}{2}\Sigma_{i=0}^{n-1}F_{xx}(t_i^n,\hat{B}_i^n)(\Delta_i^nB)^2$

$=\Sigma_{i=0}^{n-1}F_t(\hat{t}_i^n,B_{i+1}^n)\Delta_i^nt+\frac{1}{2}\Sigma_{i=0}^{n-1}F_{xx}(t_i^n,B_{i+1}^n)\Delta_i^nt+\Sigma_{i=0}^{n-1}F_x(t_i^n,B_{i+1}^n)\Delta_i^nW+\frac{1}{2}\Sigma_{i=0}^{n-1}F_{xx}(t_i^n,B_{i+1}^n)[\Delta_i^nW)^2-\Delta_i^nt]+\frac{1}{2}\Sigma_{i=0}^{n-1}[F_{xx}(t_i^n,\hat{B_i^n})-F_{xx}(t_i^n,B_i^n)](\Delta_i^nW)^2$

\begin{equation}
=A_{1,n}+A_{2,n}+A_{3,n}+A_{4,n}+A_{5,n},
\end{equation}

\noindent Since $F_t, F_x, F_{xx}$ and $B_t$ are continuous and bounded functions, we have

\begin{equation}
    lim_{n\longrightarrow\infty}sup_{i=1,...,n}sup_{t\in[t_i^n,t_{i+1}^n]}|F_t(\hat{t}_i^n,B_i^n)-F_t(t,B_t)|\longrightarrow 0
    a.s.,
\end{equation}

\begin{equation}
    lim_{n\longrightarrow\infty}sup_{i=1,...,n}sup_{t\in[t_i^n,t_{i+1}^n]}|F_{xx}(t_i^n,B_i^n)-F_{xx}(t,B_t)|\longrightarrow 0
    a.s.,
\end{equation}

\begin{equation}
    lim_{n\longrightarrow\infty}sup_{i=1,...,n}|F_{xx}(t_i^n,\hat{B}_i^n)-F_{xx}(t_i^n,B_i^n)|\longrightarrow 0
    a.s.,
\end{equation}

\noindent Now, we have a look at the sum in (3.5):

\noindent1. From (3.6), (3.7) and the definition of the Riemann integral, we have

$lim_{n\longrightarrow\infty}A_{1,n}=lim_{n\longrightarrow\infty}\Sigma_{i=0}^{n-1}F_t(\hat{t}_i^n, B_{i+1}^n)\Delta_i^nt=\int_0^TF_t(t,B_t)dt$
a.s.,and

$lim_{n\longrightarrow\infty}A_{2,n}=lim_{n\longrightarrow\infty}\Sigma_{i=0}^{n-1}F_{xx}(t_i^n, B_i^n)\Delta_i^nt=\int_0^TF_{xx}(t,B_t)dt$
a.s.

\noindent2. Since $F_x\in\mathcal{M}^2$, we can get the limit

$lim_{n\longrightarrow\infty}A_{3,n}=lim_{n\longrightarrow\infty}\Sigma_{i=0}^{n-1}F_x(t_i^n, B_i^n)\Delta_i^nB=\int_0^TF_x(t,B_t)dB_t$

according to Theorem 3.1.6.

\noindent3. For the term $A_{4,n}^2$, we have

\noindent$E(A_{4,n}^2)=E(\Sigma_{i=0}^{n-1}F_{xx}(t_i^n,B_i^n)[(\Delta_i^nB)^2-\Delta_i^nt])^2$

$=\Sigma_{i=0}^{n-1}E|F_{xx}(t_i^n,B_i^n)[(\Delta_i^nB)^2-\Delta_i^nt]|^2$
(expectation of the cross term is 0)

$=\Sigma_{i=0}^{n-1}E|F_{xx}(t_i^n,B_i^n)|^2E|\Delta_i^nB)^2-\Delta_i^nt|^2$
(the increment is independent)

$\leq C^2\Sigma_{i=0}^{n-1}E|\Delta_i^nB)^2-\Delta_i^nt|^2$
($F_{xx}$ is bounded by C)

$=2C^2\Sigma_{i=0}^{n-1}(\Delta_i^nt)^2=2C^@\Sigma_{i=0}^{n-1}\frac{T^2}{^2n}=2C^2\frac{T^2}{n}\longrightarrow 0$
as
$n\longrightarrow\infty$

\noindent4. Note that $\Sigma_{i=0}^{n-1}(\Delta_i^nB)^2\longrightarrow t$ in $\mathcal{L}^2$ and thus in probability since the left quantity is the quadratic variation of Brownian motion. Together with the continuity result (3.7), we have the following convergence(in probability):

$|A_{5,n}|=|\Sigma_{i=0}^{n-1}[F_{xx}(t_i^n,\hat{B_i^n})-F_{xx}(t_i^n,B_i^n)](\Delta_i^nW)^2|$

$\leq sup_i|F_{xx}(t_i^n,\hat{B_i^n})-F_{xx}(t_i^n,B_i^n)|\Sigma_{i=0}^{n-1}(\Delta_i^nW)^2\longrightarrow^p0$

\noindent Note that the convergence of $A_i^n$, i=1,...,5 involves different modes: $A_{1,n}$ and $A_{2,n}$ converge almost surely, $A_{3,n}$ and $A_{4,n}$ converge in $\mathcal{L}^2$, and $A_{5,n}$ converges in probability. To combine the results, note that convergence in $\mathcal{L}^2$ implies convergence in probability. Thus all $A_{3,n},A_{4,n}$ and $A_{5,n}$ converge in probability.  Note also that there is a subsequence $\{n_k\}_{k=1,2,...}$ such that $\{A_{3,n_k}\}_{k=1,2,...}$ converge a.s. Along this subsequence, we can ﬁnd a further subsequence $n_{k_l}$ such that $A_{4,n_{k_l}}$ converges a.s., and so forth. Finally, all $A_{j,n}$, j = 1,...,5 converge a.s. with respect to some subsequence $m_1 < m_2$ <..., say. Then

$F(T,B_t)-F(0,B_0)$

=$lim_{k\longrightarrow\infty}\{\Sigma_{i=0}^{n-1}F_t(\hat{t}_i^n,B_{i+1}^n)\Delta_i^nt$

$+\frac{1}{2}\Sigma_{i=0}^{n-1}F_{xx}(t_i^n,B_{i+1}^n)\Delta_i^nt+\Sigma_{i=0}^{n-1}F_x(t_i^n,B_{i+1}^n)\Delta_i^nB$

$+\frac{1}{2}\Sigma_{i=0}^{n-1}F_{xx}(t_i^n,B_{i+1}^n)[\Delta_i^nB)^2-\Delta_i^nt]+\frac{1}{2}\Sigma_{i=0}^{n-1}[F_{xx}(t_i^n,\hat{B_i^n})$

$-F_{xx}(t_i^n,B_i^n)](\Delta_i^nB)^2\}$
    
=$\int_0^T[F_t(t,B_t)+\frac{1}{2}F_{xx}(t,B_t)]dt+\int_0^TF_x(t,B_t)dB_t,$\\
\\
\noindent\textbf{Example 3.2.1}: F(t,x)=$x^3$, the partial derivatives are $F_t(t,x)=0, F_x(t,x)=3x^2$ and $F_{xx}(t,x)=6x$. According to It$\hat{o}$'s formula, we have $dB_t^3=3B_tdt+3B_t^2dB_t$ provided $3B_t^2\in\mathcal{M^2}$(which has been proved).\\
\\
\noindent Looking carefully in to the proof above, we can find that It$\hat{o}$'s Lemma also holds for F(t,$X_t$) where $X_t$ is a process with quadratic variation [$X_t$] satisfying d[$X_t$]=g(t)dt, for some g(t)$\in\mathcal{M}^2$. The following theorem gives It$\hat{o}$'s Lemma in general case.

\noindent\textbf{Theorem 3.2.2(It$\hat{o}$'s formula in general case)}: Let $X_t$ be a stochastic process with quadratic variation $[X]_t$ satisfying d[$X_t$]=g(t)dt where g(t)$\in\mathcal{M}^2$. Suppose that $F(t,x), F_t(t,x), F_x(t,x) and F_{xx}(t,x)$ are continuous for all t$\geq$0 and x$\in\mathbf{R}$. Also the process $g_tF_x(t,X_t)\in\mathcal{M}^2$. Then F(t,$X_t$) can be expressed as

\begin{equation}
    dF(T,X_t)=F_t(t,X_t)dt+F_x(t,X_t)dX_t+\frac{1}{2}F_{xx}(t,X_t)d[X]_t
\end{equation}

\noindent\textbf{Example 3.2.2}: For a process $X_t$ which satisfying $dX_t=a_tdt+b_tdB_t$, where $a_t\in\mathcal{L}^2$ and $b_t\in\mathcal{M}^2$, and $d[X]_t=(dX_t)^2=(a_tdt+b_tdB_t)^2=a_t^2(dt)^2+2a_tb_tdtdB_t+b_t^2(dB_t)^2$. Omitting the term smaller than dt, we have $d[X]_t=b_t^2dt$. Hence (3.9) reduces to

\noindent$dF(t,X_t)=F_t(t,X)dt+F_x(t,X_t)dX_t+\frac{1}{2}F_{xx}(t,X_t)d[X]_t$

$=F_t(t,X)dt+F_x(t,X_t)(a_tdt+b_tdB_t)+\frac{1}{2}F_{xx}(t,X_t)b_t^2dt$

$=(F_t(t,X)+a_tF_x(t,X_t)+\frac{1}{2}F_{xx}(t,X_t)b_t^2)dt+F_x(t,X_t)b_tdB_t$

\section{Reflecting Stochastic Differential Equation}

\subsection{Existence of solution for Skorohod Equation}

\noindent As is introduced in section two, Skorohod equation describes a reflecting Brownian motion X(t,$\omega$) on D=[0,$\infty$) which satisfies that

\begin{equation}
    X=B+\phi,
\end{equation}

\noindent where B is a standard Brownian motion and $\phi$ is a continuous stochastic process increasing only when X(t)=0. However, there is a unique solution for (4,1) not only when B is a Brownain motion but also when it is a continuous function with B(0)$\in\overline{D}$.

\noindent We firstly consider the case in which (4,1) is a multi-dimensional equation and D is a convex domain.

\noindent An $\mathbf{R}^d$-valued function $\phi(t)=(\phi^1(t),...,\phi^d(t))$ defined on $\mathbf{R}_+=[0,\infty)$ is said to be of bounded variation for simplicity if all $\phi_i(t)$ are of bounded variation on each finite t-interval. For a right continuous function $\phi(t)$ with $\phi(0)$=0, we define

\noindent$|\phi|$(t)= the total variation of $\phi$ on [0,t]

$=sup\Sigma_k|\phi(t_k)-\phi(t_{k-1})|,$

\noindent where $0=t_0<t_1<...<t_n=t$ is a partition. $\phi(t)$ can be expressed as

\begin{equation}
    \phi(t)=\int_0^t\Vec{n}(s)d|\phi|(s)=\int_{[0,t]}\Vec{n}(s)|\phi|(ds)
\end{equation}

\noindent where $\Vec{n}(t)$ is a unit vector valued function and is uniquely determined almost everywhere with respect to the measure d$|\phi|$.

\noindent We introduce some definitions as preparation

D: a convex domain in $\mathbf{R}^d$;

$\overline{D}$: closure of D;

$\mathcal{H}_x(D)$: the set of all supporting hyperplanes of D at x for x$\in\partial{D}$;

(Inward) normal vector at x$\in\partial{D}$: inward unit normal vector perpendicular to some H$\in\mathcal{H}_x(D)$;

$\mathcal{N}_x(D)$: the set of all inward normal vectors at x$\in\partial{D}$;

$\textbf{C}(\mathbf{R}_+,\mathbf{R}^b)$: the space of $\mathbf{R}^d$-valued continuous functions on $\mathbf{R}$;

$\textbf{D}(\mathbf{R}_+,\mathbf{R}^d)$: the space of $\mathbf{R}^d$-valued right continuous functions on $\mathbf{R}_+$ with left limits.

\noindent Given a function $X\in\textbf{D}(\mathbf{R_+,\overline{D}})$, a function $\phi$ is said to be associated with X if the following three conditions hold:

\noindent(i) $\phi$ is a function in $\textbf{D}(\mathbf{R}_+,\mathbf{R}^d)$ with bounded variation and $\phi(0)=0$.

\noindent(ii) The set \{$t\in\mathbf{R}_+: X(t)\in D$\} has d$|\phi|$-measure zero.

\noindent(iii) The function $\Vec{n}$(t) in (4,2) is a normal vector at X(t) for almost all t with respect to the measure d$|\phi|$. We can also use the following one as a substitute:
for any $\eta\in\textbf{C}(\mathbf{R}_+,\overline{D})$, ($\eta(t)-X(t), \phi(dt))\geq$0.

\noindent Then our problem can be expressed as follows:

\noindent Given w$\in\textbf{D}(\mathbf{R}_+,\mathbf{R}^d)$ with w(0)$\in\overline{D}$, find a solution X of

\begin{equation}
    X=w+\phi,
\end{equation}

\noindent and it is always assumed that X$\in\textbf{D}(\mathbf{R}_+,\overline{D})$ and $\phi$ is associated with X.\\
\\
\noindent In the general multi-dimensional case, the existence of a solution of (4.3) is not trivial. However when w is a step function, we can easily find the solution. For a given point x$\in\mathbf{R}^d-\overline{D}$, we denote the (unique) point on $\partial{D}$ which gives the minimum distance between x and $\overline{D}$ by $[x]_\partial$.

\noindent\textbf{Lemma 4.1.1} \cite{tanaka2002stochastic}: If w is a step function with w(0)$\in\overline{D}$, there exists a solution of (4.3).

\noindent\textbf{Proof}: To prove the existence, we can try to construct a X(t) satisfying the conditions. Luckily in the case, it is not hard to do that. Let $T_1=inf\{t>0:w(t)\notin\overline{D}$ and then define X to be X(t)=$w$(t) for $0\leq t<T_1$ and X$(T_1)=[w(T_1)]_\partial$, in this case $\phi(t)=0$ for t<$T_1$ and $\phi(T_1)=[w(T_1)]_\partial-w(T_1)$. It solves (4.3) for $0\leq t\leq T_1$. Then we consider t$>T_1$. Suppose that the solution of X(t) has been found for $0\leq t\leq T_{n-1}$. Let $T_n=inf{t>T_{n-1}:w(t)+\phi(T_{n-1})\notin\overline{D}}$,

\[X(t)=
\begin{cases}
\ w(t)+\phi(T_{n-1}), & for T_{n-1}<t<T_n;\\
\ [w(T_n)+\phi(T_{n-1})]_\partial, & for t=T_n .
\end{cases}\]

\noindent By repeating this step we will get the solution of (4.3) since $T_n\uparrow\infty$ as $n\uparrow\infty$.

\noindent\textbf{Lemma 4.1.2} \cite{tanaka2002stochastic}:(i) Let $w,\Tilde{w}\in\textbf{D}(\mathbf{R}_+,\mathbf{R}^d)$ with w(0),$\Tilde{w}(0)\in\overline{D}$ and X,$\Tilde{X}$ be the solution of X=$w+\phi$, $\Tilde{X}=\Tilde{w}+\Tilde{\phi}$ respectively. Then we have

$|X(t)-\Tilde{X}(t)|^2\leq|w(t)-\Tilde{w}(t)|^2+2\int_0^t(w(t)-\Tilde{w}(t)-w(s)+\Tilde{w}(s))(\phi(ds)-\Tilde{\phi}(ds))$.

\noindent(ii) If X is a solution of (4.3), then

$|X(t)-X(s)|^2\leq|w(t)-w(s)|^2+2\int_{(s,t]}(w(t)-w(\tau)\phi(d\tau), 0\leq s\leq t$.

\noindent Now we can prove the uniqueness of the solution.

\noindent\textbf{Lemma 4.1.3} \cite{tanaka2002stochastic}: (4.3) has at most one solution.

\noindent\textbf{Poof}: Suppose X and $\Tilde{X}$ are both solutions of (4.3) and set $w=\Tilde{w}$. By (i) of lemma(4.2) we have $|X(t)-\Tilde{X}(t)|\leq|w(t)-\Tilde{w}(t)|^2+2\int_0^t(w(t)-\Tilde{w}(t)-w(s)+\Tilde{w}(s))(\phi(ds)-\Tilde{\phi}(ds))=0$.

\noindent\textbf{Lemma 4.1.4} \cite{tanaka2002stochastic}: If $w$ is continuous, then the solution of (4.4) is also continuous.

\noindent Use the inequation in (ii) of lemma(4.2), it is easy to prove the lemma.

\noindent\textbf{Lemma 4.1.5} \cite{tanaka2002stochastic}: Let $\{w_n\}_{n\geq1}$ be a sequence in $\textbf{D}(\mathbf{R}_+,\mathbf{R}^d)$ such that for each n the equation $X_n=w_n+\phi_n$ has a solution for $0\leq t\leq T$, T being a positive constant. If $w_n$ converges uniformly on [0,T] to some $w\in\textbf{C}(\mathbf{R}_+,\mathbf{R}^d)$ as $n\longrightarrow\infty$ and if $\{|\phi_n|(T)\}_{n\geq1}$ is bounded, then $X_n$ converges uniformly on [0,T] as $n\longrightarrow\infty$ to the solution X=$w+\phi$ for $0\leq t\leq T$.

\noindent\textbf{Proof}: Assume $|\phi_n|(t)\leq K$ for each n. Then by (i) of lemma(4.2), we have

\noindent$|X_n(t)-X_m(t)|^2\leq|w_n(t)-w_m(t)|^2$

$+2\int_0^t(w_n(t)-w_m(t)-w_n(s)+w_m(s))(\phi_n(ds)-\phi_m(ds))$

$\leq|w_n(t)-w_m(t)|^2+8K {sup}_{0\leq s\leq t}|w_n(s)-w_m(s)|$. 

\noindent Since the sequence $\{w_n\}_{n\geq1}$ is uniformly convergent to $w$, $X_n=w_n+\phi_n$ also converges uniformly to X=$w+\phi$. Then we prove that X is the solution of (4.3). To prove this is just to prove $\phi$ is associated with X. The condition (i) is satisfied obviously since $|\phi_n|$(T) is bounded. Condition (ii) is also trivial. Then we try to verify condition (iii). Let $\eta\in\textbf{C}(\mathbf{R}_+,\overline{D})$ and notice that for $0\leq t_1<t_2\leq T$

$|\int_{t_1}^{t_2}(\eta(t)-X_n(t))\phi_n(dt)-\int_{t_1}^{t_2}(\eta(t)-X(t))\phi(dt)\leq$

$|\int_{t_1}^{t_2}(X(t)-X_n(t))\phi_n(dt)|+|\int_{t_1}^{t_2}(\eta(t)-X(t))(\phi_n(dt)-\phi(dt))|$.

\noindent The first is dominated by K${sup}_{t_1\leq t\leq t_2}|X(t)-X_n(t)|$ and hence tends to 0 as n$\longrightarrow\infty$; the second term also does as can be seen by approximating the integral by the Riemann sum. Therefore

$\int_{t_1}^{t_2}(\eta(t)-X(t))\phi(dt)=lim_{n\longrightarrow\infty}\int_{t_1}^{t_2}(\eta(t)-X_n(t))\phi_n(dt)\geq0$

\noindent and the proof finished.\\
\\
\noindent Now we move on to the problem of existence of solution for (4.3) under the assumption that $w\in\textbf{C}(\mathbf{R}_+,\mathbf{R}^d)$.

\noindent Here, we only introduce two special conditions where the solution exists.

\noindent\textbf{Condition A}: There is a unit vector $\Vec{e}$ and a constant c\textgreater0 such that ($\Vec{e},\Vec{n})\geq$c for any $\Vec{n}\in\cup_{y\in\partial{d}}\mathcal{N}_y$(\textbf{D}.

\noindent\textbf{Lemma 4.1.6} \cite{tanaka2002stochastic}: Assume that D satisfies the condition A. Then there exists a solution X of (4.3) for any $w\in\textbf{C}(\mathbf{R}_+,\mathbf{R}^d)$, and for $0\leq s<$t

$|X(t)-X(s)|\leq K\Delta_{s,t}$,

$|\phi|(t)-|\phi|(s)\leq K'\Delta_(s,t)$,

\noindent where K and K' are constant depending only on the constant c in the condition(A); $\Delta_{s,t}={sup}_{s\leq t_1<t_2\leq t}|w(t_2)-w(t_1)|$.\\
\\
\noindent\textbf{Condition B}: There exist X\textgreater0 and $\delta>$0 such that for any $x\in\partial{D}$ we can find an open ball $B_\epsilon(x_0)={y\in\mathbf{R}^d:|y-x_0|<\epsilon}$ satisfying $B_\epsilon(x_0)\subset D$ and $|x-x_0|\leq\delta$.

\noindent Condition B is always satisfied if D is bounded or if d=2.

\noindent\textbf{Lemma 4.1.7} \cite{tanaka2002stochastic}: Assume that D satisfies the condition B. Then there exists a unique solution of (4.3) if $w\in\textbf{C}(\mathbf{R}_+,\mathbf{R}^d)$, and the solution depends continuously on $w$ with respect to the compact uniform topology.

\noindent\textbf{Theorem 4.1.8} \cite{tanaka2002stochastic}: Let D be a general convex domain and $\{w_n\}_{n\geq1}$ be a sequence in $\textbf{C}(\mathbf{R}_+,\mathbf{R}^d)$ such that $X_n=w_n+\phi_n$ has a solution for each n. Assume that $w_n$ and $X_n$ converge to $w$ and X uniformly on compacts as $n\longrightarrow\infty$, respectively. Then X is a solution of (4.3).

\noindent\textbf{Proof}: For any constant T\textgreater0, there is a constant N such that

\[{sup}_n{max}_{0\leq t\leq T}|X_n(t)|<N\]

\noindent For such N both $X_n$ and X are the solutions of (4.3) when $0\leq t\leq T$. Then we construct a domain $D_N=D\cap\{|x|<N\}$ and it satisfies condition B. Hence by lemma 4.7 X is the solution of (4.3) for $D_N$ and so for D.

\subsection{Stochastic version of Skorohod equation}

\noindent We aim to verify the existence of solution of (4.3) without satisfying condition B.

\noindent Let $(\Omega,\mathcal{F},P)$ be a complete probability space with an increasing family $\{\mathcal{F}\}_{t\geq 0}$ of sub-$\sigma$-fields of $\mathcal{F}$ where $\mathcal{F}_t=\cup_{\epsilon>0}\mathcal{F}_{t+\epsilon}$. Let D be a convex domain. Then we will extend the theorems and lemmas before to stochastic case.

\noindent\textbf{Theorem 4.2.1} \cite{tanaka2002stochastic}: Let {M(t)} be an $\mathbf{R}^d$-valued process with M(0)$\in\overline{D}$ such that each component is a continuous local $\mathcal{F}_t$-martingale and {A(t)} be an $\mathbf{R}^d$ valued, continuous and $\mathcal{F}_t$-adapted process of bounded variation with A(0)=0. then there exists a unique $\mathcal{F}_t$-adapted solution \{X(t)\} of

\begin{equation}
    X(t)=M(t)+A(t)+\Phi(t).
\end{equation}

\noindent Moreover, for $f\in\mathbf{C}^2(\mathbf{R}$ with $f'\geq0$ on $\mathbf{R}_+$ and $0\leq s\leq t$ we have

\begin{equation}
    f(|X(t)-X(s)|^2)
\end{equation}

\[\leq f(0)+2\int_s^tf'\Sigma_i(X^i(\tau)-X^i(s))(M^i(d\tau)+A^i(d\tau))\]

\[+2\int_s^tf''\Sigma_{i,j}(X^i(\tau)-X^i(s))(X^j(\tau)-X^j(\tau))d[M^i,M^j]\]

\[+\int_s^tf'\Sigma_id[M^i,M^j],\]

\noindent where f' and f'' are evaluated at $|X(\tau)-X(s)|^2$ and $[M^i,M^j]$ denotes the quadratic variation process.\\
\\
\noindent By a solution of (4.4), we mean a $\overline{D}$-valued process \{X(t)\} which satisfies (4.4) almost surely, under the condition that almost all sample paths of $\{\Phi(t)\}$ are associated with those of \{X(t)\}.

\subsection{Stochastic differential equation with reflection}

\noindent Let D be a convex domain in $\mathbf{R}^d$ and $\{\Omega,\mathcal{F},P;\mathcal{F}_t\}$ satisfy the same condition as in the last subsection. B(t)=($B^1(t),...,B^r(t)$) with B(0)=0 is an $\mathcal{F}_t$-adapted r-dimensional Brownian motion and for $0\leq s\leq t,\varepsilon\in\mathbf{R}^d$

\[E[e^{i(\varepsilon,B(t)-B(s))}|\mathcal{F}_s]=e^\frac{-(t-s)|\varepsilon|^2}{2}, a.s.\]

\noindent Let $\sigma(t,x)=\{\sigma_k^i(t,x)\}$ be an $\mathbf{R}^d\bigotimes\mathbf{R}^r$-valued function and $b(t,x)=\{b^i(t,x)\}$ be an $\mathbf{R}^d$-valued function, both being defined on $\mathbf{R}_+\times\overline{D}$. Then we consider the stochastic differential equation with reflection

\begin{equation}
    dX=\sigma(t,X)dB+b(t,X)st+d\Phi, X(0)=x,
\end{equation}

\noindent where x=$(x^1,...,x^d)\in\overline{D}$. The aim is to find an $\mathcal{F}_t$-adapted $\overline{D}$-process \{X(t)\} under the condition that $\{\Phi(t)\}$ is an associated process of \{X(t)\}. $\sigma(t,x)$ and b(t,x) are always assumed to be Borel measurable in (t,x).

\noindent\textbf{Theorem 4.3.1} \cite{tanaka2002stochastic}: If there exists a constant K\textgreater0 such that

\begin{equation}
    ||\sigma(t,x)-\sigma(t,x)||\leq K|x-y|, |b(t,x)-b(t,y)|\leq K|x-y|,
\end{equation}

\begin{equation}
    ||\sigma(t,x)||\leq K(1+|x|^2)^\frac{1}{2}, ||b(t,x)||\leq K(1+|x|^2)\frac{1}{2}
\end{equation}

\noindent then there exists a (pathwise) unique $\mathcal{F}_t$-adapted solution of (4.6) for any x$\in\overline{D}$.

\noindent Before proving theorem 4.3.1, we introduce an inequation first.

\noindent\textbf{Lemma 4.3.2} \cite{tanaka2002stochastic}: Replace $w,\Tilde{w}$ in (i) of lemma 4.1.2 by $w+a,\Tilde{w}+\Tilde{a}$, respectively, where a and $\Tilde{a}$ are $\mathbf{R}^d$-valued right continuous functions of bounded variation with a(0)=$\Tilde{a}$(0)=0, then

\[|X(t)-\Tilde{X}(t)|^2\leq|w(t)-\Tilde{w}(t)|^2+2\int_0^t(X(s)-\Tilde{X}(s))(a(ds)-\Tilde{a}(ds))\]

\[+2\int_0^t(w(t)-\Tilde{w}(t)-w(s)+\Tilde{w}(s))(a(ds)-\Tilde{a}(ds)).\]

\noindent By a similar replacement of $w$ in (ii) of lemma 4.1.2 by $w+a$, we have the following inequation

\[|X(t)-X(s)|^2\leq|w(t)-w(s)|^2+2\int_0^t(X(\tau)-X(s))a(d\tau)\]

\[+2\int_{(s,t]}(w(t)-w(\tau))(a(d\tau)+\phi(d\tau)).\]

\noindent\textbf{Theorem 4.3.3} \cite{tanaka2002stochastic}: If $\sigma$(t,x) and b(t,x) are bounded continuous on $\mathbf{R}_+\times\overline{D}$, then on some probability space ($\Omega,\mathcal{F},P$) we can find an r-dimensional Brownain motion \{B(t)\} in such a way that (4.6) has a solution.



\bibliographystyle{elsarticle-harv} 
\bibliography{elsarticle-template-harv}

\begin{thebibliography}{17}
\expandafter\ifx\csname natexlab\endcsname\relax\def\natexlab#1{#1}\fi
\providecommand{\url}[1]{\texttt{#1}}
\providecommand{\href}[2]{#2}
\providecommand{\path}[1]{#1}
\providecommand{\DOIprefix}{doi:}
\providecommand{\ArXivprefix}{arXiv:}
\providecommand{\URLprefix}{URL: }
\providecommand{\Pubmedprefix}{pmid:}
\providecommand{\doi}[1]{\href{http://dx.doi.org/#1}{\path{#1}}}
\providecommand{\Pubmed}[1]{\href{pmid:#1}{\path{#1}}}
\providecommand{\bibinfo}[2]{#2}
\ifx\xfnm\relax \def\xfnm[#1]{\unskip,\space#1}\fi
\bibitem[{Asmussen et~al.(1995)Asmussen, Glynn and
  Pitman}]{asmussen1995discretization}
\bibinfo{author}{Asmussen, S.}, \bibinfo{author}{Glynn, P.},
  \bibinfo{author}{Pitman, J.}, \bibinfo{year}{1995}.
\newblock \bibinfo{title}{Discretization error in simulation of one-dimensional
  reflecting brownian motion}.
\newblock \bibinfo{journal}{The Annals of Applied Probability} ,
  \bibinfo{pages}{875--896}.
\bibitem[{Berestycki et~al.(2014)Berestycki, Berestycki and
  Schweinsberg}]{berestycki2014critical}
\bibinfo{author}{Berestycki, J.}, \bibinfo{author}{Berestycki, N.},
  \bibinfo{author}{Schweinsberg, J.}, \bibinfo{year}{2014}.
\newblock \bibinfo{title}{Critical branching brownian motion with absorption:
  survival probability}.
\newblock \bibinfo{journal}{Probability Theory and Related Fields}
  \bibinfo{volume}{160}, \bibinfo{pages}{489--520}.
\bibitem[{Calin(2012)}]{calin2012introduction}
\bibinfo{author}{Calin, O.}, \bibinfo{year}{2012}.
\newblock \bibinfo{title}{An introduction to stochastic calculus with
  applications to finance}.
\newblock \bibinfo{journal}{Ann Arbor} .
\bibitem[{D’Auria and Kella(2012)}]{d2012markov}
\bibinfo{author}{D’Auria, B.}, \bibinfo{author}{Kella, O.},
  \bibinfo{year}{2012}.
\newblock \bibinfo{title}{Markov modulation of a two-sided reflected brownian
  motion with application to fluid queues}.
\newblock \bibinfo{journal}{Stochastic Processes and their Applications}
  \bibinfo{volume}{122}, \bibinfo{pages}{1566--1581}.
\bibitem[{Grebenkov(2019)}]{grebenkov2019probability}
\bibinfo{author}{Grebenkov, D.S.}, \bibinfo{year}{2019}.
\newblock \bibinfo{title}{Probability distribution of the boundary local time
  of reflected brownian motion in euclidean domains}.
\newblock \bibinfo{journal}{Physical Review E} \bibinfo{volume}{100},
  \bibinfo{pages}{062110}.
\bibitem[{Guo and Li(2019)}]{guo2019novel}
\bibinfo{author}{Guo, X.}, \bibinfo{author}{Li, J.}, \bibinfo{year}{2019}.
\newblock \bibinfo{title}{A novel twitter sentiment analysis model with
  baseline correlation for financial market prediction with improved
  efficiency}, in: \bibinfo{booktitle}{2019 Sixth International Conference on
  Social Networks Analysis, Management and Security (SNAMS)},
  \bibinfo{organization}{IEEE}. pp. \bibinfo{pages}{472--477}.
\bibitem[{Ikeda and Watanabe(2014)}]{ikeda2014stochastic}
\bibinfo{author}{Ikeda, N.}, \bibinfo{author}{Watanabe, S.},
  \bibinfo{year}{2014}.
\newblock \bibinfo{title}{Stochastic differential equations and diffusion
  processes}.
\newblock \bibinfo{publisher}{Elsevier}.
\bibitem[{Kanagawa(2009)}]{kanagawa2009numerical}
\bibinfo{author}{Kanagawa, S.}, \bibinfo{year}{2009}.
\newblock \bibinfo{title}{Numerical analysis of reflecting brownian motion and
  a new model of semi-reflecting brownian motion with some domains}.
\newblock \bibinfo{journal}{Communications in Applied Analysis}
  \bibinfo{volume}{13}, \bibinfo{pages}{231}.
\bibitem[{Lang(1995)}]{lang1995effective}
\bibinfo{author}{Lang, R.}, \bibinfo{year}{1995}.
\newblock \bibinfo{title}{Effective conductivity and skew brownian motion}.
\newblock \bibinfo{journal}{Journal of statistical physics}
  \bibinfo{volume}{80}, \bibinfo{pages}{125--146}.
\bibitem[{Li(2020)}]{li2020low}
\bibinfo{author}{Li, J.}, \bibinfo{year}{2020}.
\newblock \bibinfo{title}{Low-loss tunable dielectrics for millimeter-wave
  phase shifter: from material modelling to device prototyping}, in:
  \bibinfo{booktitle}{IOP Conference Series: Materials Science and
  Engineering}, \bibinfo{organization}{IOP Publishing}. p.
  \bibinfo{pages}{012057}.
\bibitem[{Li and Guo(2020)}]{li2020global}
\bibinfo{author}{Li, J.}, \bibinfo{author}{Guo, X.}, \bibinfo{year}{2020}.
\newblock \bibinfo{title}{Global deployment mappings and challenges of
  contact-tracing apps for covid-19}.
\newblock \bibinfo{journal}{Available at SSRN 3609516} .
\bibitem[{Malliaris(1983)}]{malliaris1983ito}
\bibinfo{author}{Malliaris, A.}, \bibinfo{year}{1983}.
\newblock \bibinfo{title}{It{\^o}’s calculus in financial decision making}.
\newblock \bibinfo{journal}{SIAM review} \bibinfo{volume}{25},
  \bibinfo{pages}{481--496}.
\bibitem[{Malsagov and Mandjes(2019)}]{malsagov2019approximations}
\bibinfo{author}{Malsagov, A.}, \bibinfo{author}{Mandjes, M.},
  \bibinfo{year}{2019}.
\newblock \bibinfo{title}{Approximations for reflected fractional brownian
  motion}.
\newblock \bibinfo{journal}{Physical Review E} \bibinfo{volume}{100},
  \bibinfo{pages}{032120}.
\bibitem[{Reiman(1984)}]{reiman1984open}
\bibinfo{author}{Reiman, M.I.}, \bibinfo{year}{1984}.
\newblock \bibinfo{title}{Open queueing networks in heavy traffic}.
\newblock \bibinfo{journal}{Mathematics of operations research}
  \bibinfo{volume}{9}, \bibinfo{pages}{441--458}.
\bibitem[{Saisho(1987)}]{saisho1987stochastic}
\bibinfo{author}{Saisho, Y.}, \bibinfo{year}{1987}.
\newblock \bibinfo{title}{Stochastic differential equations for
  multi-dimensional domain with reflecting boundary}.
\newblock \bibinfo{journal}{Probability Theory and Related Fields}
  \bibinfo{volume}{74}, \bibinfo{pages}{455--477}.
\bibitem[{Tanaka(2002)}]{tanaka2002stochastic}
\bibinfo{author}{Tanaka, H.}, \bibinfo{year}{2002}.
\newblock \bibinfo{title}{Stochastic differential equations with reflecting
  boundary condition in convex regions}, in: \bibinfo{booktitle}{Stochastic
  Processes: Selected Papers of Hiroshi Tanaka}. \bibinfo{publisher}{World
  Scientific}, pp. \bibinfo{pages}{157--171}.
\bibitem[{Zhang(2020)}]{zhang2020value}
\bibinfo{author}{Zhang, Y.}, \bibinfo{year}{2020}.
\newblock \bibinfo{title}{The value of monte carlo model-based variance
  reduction technology in the pricing of financial derivatives}.
\newblock \bibinfo{journal}{PloS one} \bibinfo{volume}{15},
  \bibinfo{pages}{e0229737}.

\end{thebibliography}





\end{document}